\documentclass [12pt]{amsart}
\newif\ifpdf
\ifx\pdfoutput\undefined
	\pdffalse
\else
	\pdfoutput=1
	\pdftrue
\fi
\usepackage{enumerate}
\usepackage{amssymb}
\ifpdf
	\usepackage[pdftex]{graphicx}
	\DeclareGraphicsExtensions{.pdf}
	\usepackage{hyperref}
\else
	\usepackage{graphicx}
	\DeclareGraphicsExtensions{.eps}
	\usepackage[hypertex]{hyperref}
\fi
\usepackage{amscd}
\newcommand{\FF}{\mathbb F}

\newcommand{\ZZ}{\mathbb Z}
\newcommand{\PP}{\mathbb P}

\newcommand{\CC}{\mathbb C}
\newcommand{\NN}{\mathbb N}

\newcommand{\TT}{\mathbb T}
\newcommand{\bbA}{\mathbb A}

\newcommand{\NS}{\mathop {\rm NS}\nolimits}

\newcommand{\Pic}{\mathop {\rm Pic}\nolimits}
\newcommand{\Supp}{\mathop {\rm Supp}\nolimits}
\newcommand{\Spec}{\mathop {\rm Spec}\nolimits}
\newcommand\pior{\pi_1^{\text{\rm orb}}}

\newcommand\surj{\twoheadrightarrow}

\newcommand{\scU}{\mathcal U}
\newcommand{\scV}{\mathcal V}

\newcommand{\scP}{\mathcal P}

\newcommand{\Sing}{\mathop {\rm Sing}\nolimits}

\newtheorem{thm0}{Theorem}
\newtheorem{cor0}[thm0]{Corollary}
\newtheorem{thm}{Theorem}[section]
\newtheorem{cor}[thm]{Corollary}
\newtheorem{prop}[thm]{Proposition}
\newtheorem{lem}[thm]{Lemma}

\newtheorem*{claim}{Claim}
\theoremstyle{definition}
\newtheorem{defin}[thm]{Definition}
\newtheorem{exmple}[thm]{Example}
\newtheorem{exmples}[thm]{Examples}
\theoremstyle{remark}
\newtheorem{rem}[thm]{Remark}
\newtheorem{rem0}[thm0]{Remark}

\title{Pencils and Infinite Dihedral covers of $\PP^2$}

\author[E. Artal]{Enrique Artal Bartolo}

\author[J.I. Cogolludo]{Jos\'e Ignacio Cogolludo}
\address{Departamento de Matem\'aticas\\
Universidad de Zaragoza\\
Campus Plaza San Francisco s/n\\
E-50009 Zaragoza SPAIN}
\email{artal@unizar.es,jicogo@unizar.es}

\author[H. Tokunaga]{Hiro-o Tokunaga}
\address{Department of Mathematics \\
Tokyo Metropolitan University \\
Minamiohsawa Hachoji  \\
192-0357  Tokyo JAPAN}
\email{tokunaga@comp.metro-u.ac.jp}
\thanks{First author and second author are partially
supported by 
BFM2001-1488-C02-02}

\begin{document}
\begin{abstract}
In this work we study the connection between the existence of finite dihedral covers of 
the projective plane ramified along an algebraic curve $C$, infinite dihedral covers, and
pencils of curves containing $C$.
\end{abstract}
\maketitle

\section*{Introduction}

Let us consider a reduced plane curve $C\subset\PP^2$. The third author has extensively studied algebraic conditions for the existence of dihedral covers of $\PP^2$ ramified along $C$. In this paper, $C$ will be supposed to have two irreducible components $C_1$ and $C_2$
with the purpose to study the existence of $D_{2n}$-covers of $\PP^2$ branched at $2C_1 + n C_2$, 
for $n$ odd (see the comments before Theorem~\ref{thm:main} for the notations). 
Such covers are related to epimorphisms $\pi_1(\PP^2\setminus C)\to D_{2 n}$
sending a meridian of $C_1$ (resp. $C_2$) to a conjugate of
$\sigma$ (resp. $\tau$), see subsection~\ref{ss-dih-cov}.
Our goal is to derive the existence of $(\ZZ/2*\ZZ/2)$-covers that factorize through
such finite dihedral covers. This will be related to the existence of pencils of curves
containing $C$ and the existence of infinite dihedral covers of $\PP^2$.
We will impose some restrictions on the curves $C$; some of them are necessary conditions
for the existence of the above $D_{2n}$-covers and others will be set for the sake of
simplicity.

\begin{enumerate}[(i)]
\item $\deg C_1$ is even: this is a necessary condition for the existence
of the intermediate double cover ramified along $C_1$, see subsection~\ref{ss-dih-cov}.
\item $C_1$ has at most simple singularities: this condition will simplify some proofs.
\item $C_2\cap \Sing(C_1) = \emptyset$.
\item For each local branch $\varphi$ of $C_2$ at $P\in C_1\cap C_2$, $(\varphi\cdot C_1)_P$ 
is even: this is also a necessary condition for the reducibility of the preimage of $C_2$ by
the double cover ramified on $C_1$, see Proposition~\ref{prop:necdi}.
\end{enumerate}

Let us introduce the general setting of this work. Let $X$ and $Y$ be normal projective
varieties. Let $\pi : X \to Y$ be a finite surjective morphism. Under these conditions, 
the rational function field $\CC(X)$ of $X$ is regarded as a field extension of $\CC(Y)$, 
the function field of $Y$. We call $X$ a {\emph{$D_{2n}$-cover}} of $Y$ if the field extension
$\CC(X)/\CC(Y)$ is Galois and its Galois group is isomorphic to the dihedral group $D_{2n}$ 
of order $2n$. 

The branch locus  of  $\pi : X \to Y$, denoted by $\Delta(X/Y)$ or
$\Delta_{\pi}$ is the subset of $Y$ given by
$$
\Delta_{\pi} := \{ y\in Y \mid \pi \text{ is not a local isomorphism at } y \}.
$$
It is well known that $\Delta_{\pi}$ is an algebraic subset of  codimension $1$ if $Y$ is
smooth, see~\cite{zariski4}. Suppose that $Y$ is smooth and let 
$\Delta_{\pi} =  B_1+\dots + B_r$ be its irreducible decomposition. 
We say $\pi : X \to Y$ is branched at $e_1B_1+ \cdots +e_rB_r$ if 
the ramification index along $B_i$ is~$e_i$.

Let us state our main results:

\begin{thm0}\label{thm:main}
If $D_{2n}$-covers of $\PP^2$ branched at $2C_1 + n C_2$ exist for  enough odd numbers
$n \in \NN$, then they exist for any $n \in \NN$. Moreover, if $F_i$ denote defining 
equations of $C_i$, $i=1,2$, then there exist homogeneous polynomials $G_1$ and $G_2$ 
such that $F_2 = G_1^2 - G_2^2F_1$.
\end{thm0}

\begin{cor0}\label{cor:main}
Under the hypothesis of Theorem~{\rm\ref{thm:main}}, there exists an epimorphism from 
$\pi_1(\PP^2\setminus (C_1\cup C_2))$ onto the infinite dihedral group $\ZZ_2*\ZZ_2$.
\end{cor0}

\begin{rem0}\label{rem-num}
It is possible to be more precise in the statement of Theorem~{\rm\ref{thm:main}} in terms
of the curves $C_1,C_2$. Consider the standard resolution of the double cover of
$\PP^2$ ramified along $C_1$. By Proposition~{\rm\ref{prop:necdi}} the preimage of $C_2$ 
under this cover decomposes as $C_2^+\cup C_2^-$ into two irreducible components. 
As shown in equation~\eqref{ast1}, divisibility properties of $C_2^+ - C_2^-$ are 
required for Theorem~{\rm\ref{thm:main}} to hold. For instance, let $\nu$ be the
self-intersection of $C_2^+ - C_2^-$ and assume that $\nu\neq 0$, then the existence of a 
single $D_{2n}$-cover of $\PP^2$ branched at $2C_1 + n C_2$, where $n^2$ does not divide 
$\nu$, is enough for Theorem~{\rm\ref{thm:main}} to hold.
\end{rem0}

%\input inf_dihed_sec1.tex
%inf_dihed_sec1.tex

\section{Preliminaries}

\subsection{Topology of a double cover of $\PP^2$}
\mbox{}

Let $B$ be a reduced plane curve of even degree $d$. Assume that singularities of 
$B$ are all simple.
Let  $\delta : Z \to \PP^2$ be a double cover branched
at $B$ and let $\mu :\tilde Z \to Z$ be the canonical resolution, see
\cite{horikawa}. 

\begin{lem}\label{lem:P2} { $\tilde Z$ is simply connected.}
\end{lem}

\begin{proof}
By using results on the simultaneous resolution (\cite{brieskorn1,brieskorn2}), 
we may assume that $B$ is smooth. In this case, $\tilde Z = Z$.
If $B$ is smooth, $\pi_1(\PP^2 \setminus B) \cong \ZZ/d\ZZ$. Hence
$\pi_1(Z \setminus {\delta}^{-1}(B)) \cong \ZZ/(d/2)\ZZ$ and it is generated by a 
meridian around ${\delta}^{-1}(B)$. In $Z$, this lasso is homotopic to zero. Hence
$\pi_1(Z) = \{1\}$.
\end{proof}

\begin{cor}\label{cor:P2} { $\Pic(\tilde Z) = \NS(\tilde Z)$,
$\Pic(\tilde Z)$ is a lattice with respect to the intersection pairing.}
\end{cor}

\subsection{Dihedral covers}\label{ss-dih-cov}
\mbox{}

To present $D_{2n}$, we use the notation
\[
D_{2n} = \langle \sigma, \tau \mid \sigma^2 = \tau^n = (\sigma\tau)^2 = 1 \rangle,
\]
and fix it throughout this article. Given a $D_{2n}$-cover $\pi : X \to Y$, we canonically
obtain the double cover $D(X/Y)$ of $Y$ by taking the $\CC(X)^{\tau}$-normalization of
$Y$, where $\CC(X)^{\tau}$ is the fixed field of $\langle \tau \rangle$. The variety $X$ is an 
$n$-cyclic cover of $D(X/Y)$ by its definition. We denote these cover morphisms by
$\beta_1(\pi) : D(X/Y) \to Y$ and $\beta_2(\pi) : X \to D(X/Y)$, respectively. 
In this context, $D_{2n}$-covers have been thoroughly considered in \cite{tokunaga1,tokunaga2}. We here summarize some of the main results that will be
needed in this paper. 
The following is a sufficient condition for the existence of $D_{2n}$-covers.

\begin{prop}\label{prop:con_di-1}
Let $Y$ be a smooth variety and let $n\geq 3$ an integer.
Let $\delta :Z \to Y$ be a smooth double cover, and let $\sigma_\delta$ 
denote the involution. Suppose that there exists a reduced divisor
$D$ on $Z$ satisfying the following conditions:

\begin{enumerate}[\rm(1)]

\item $\sigma_\delta^*D$ and  $D$ have no common components.

\item  There exists a line bundle $L$ on $Z$ such that
$D  - \sigma_\delta^*D \sim n L$, where $\sim$ means linear equivalence.
\end{enumerate}
Let us suppose also that either $n$ is odd, or $n$ is even and
$Y$ is simply connected.
Then there  exists a $D_{2n}$-cover $\pi : X \to Y$ such that
$D(X/Y) = Z$ and $\pi$ is branched at $2\Delta_\delta + n\delta(D)$.
\end{prop}

\begin{proof}  If $n$ is odd,
our statement is a special case of \cite{tokunaga2}. If $n$ is
even, then by
\cite[Remark 3.1]{tokunaga1}
and a similar argument to the proof of
\cite[Proposition 1.1]{tokunaga2}, the result follows.
\end{proof}

\begin{cor}\label{cor:con_di-1}{ Suppose that $Y$ is simply connected. If $\sigma_f^*D \sim D$, then
there exists a $D_{2n}$-cover of $Y$ branched at $2\Delta_f + nf(D)$ for any $n \ge 3$.
}
\end{cor}

As for a necessary condition for the existence  of
$D_{2n}$-covers, we have the following.

\begin{prop}[{\cite[\S2]{tokunaga1}}]
\label{prop:necdi}   
Let $\pi : X \to Y$ be a $D_{2n}$-cover such
that $D(X/Y)$ is smooth. Then there exist a (possibly empty) effective
divisor $D_1$ and a line bundle $L$ on $D(X/Y)$ satisfying
the following conditions:
\begin{enumerate}[\rm(1)]
\item $D_1$ and $\sigma^*D_1$ have no common components. 

\item $D_1 - \sigma^*D_1 \sim n L$.

\item $\Delta(X/D(X/Y)) = \Supp(D_1 + \sigma^*D_1)$.

\item The ramification index along $D_{1,j}$ is $\frac{n}{\gcd(a_j,n)}$,
where $D_1=\sum_j a_jD_{1,j}$, $(a_j>0)$ is the irreducible decomposition.

\end{enumerate}
\end{prop}

\begin{cor}\label{cor:branch}  Let $D$ be an irreducible component of $\beta_1(\pi)(\Delta_{\beta_2(\pi)})$.
Then the divisor $\beta_1(\pi)^*D$ is of the form $D' + \sigma^*D'$ for some irreducible divisor on $D(X/Y)$. In other
words, $\beta_2(\pi)$ is not branched along any irreducible divisor $D$ with $D = \sigma^*D$.
\end{cor}

%\input inf_dihed_sec2.tex
%inf_dihed_sec2.tex

\section{Certain $D_{2n}$-covers of algebraic surfaces}

Let $\Sigma_o$ be a smooth projective surface. Let $C_1$ and $C_2$
be reduced divisors on $\Sigma_o$ such that 
\begin{itemize}
\item $C_1$ has at most
simple singularities;
\item $C_2$ is irreducible;
\item $C_2 \cap \Sing(C_1) = \emptyset$;
\item there exists a double 
cover $\delta : Z \to \Sigma_o$ branched at $C_1$;
\item its canonical resolution $\mu : \tilde Z \to Z$ is simply connected.
\end{itemize} 
\medskip

\begin{prop}\label{prop:key}{If there exists a $D_{2k}$-cover
$\pi_k : S_k \to \Sigma_o$
branched at $2C_1 + kC_2$ for finitely many enough odd natural numbers
$k$ (see Remark~{\rm\ref{rem-num}}), then
there exist $D_{2n}$-covers of $\Sigma_o$ branched at $2C_1 + nC_2$ for 
any integer $n \ge 3$}
\end{prop}

\begin{proof}
By our assumption, $D(S_k/\Sigma_o) = Z$ and $\beta_1(\pi_k) = \delta$. Let 
$$
\begin{CD}
Z & @<{\mu}<< & \tilde Z \\
@V{\delta}VV&  & @VV{\tilde\delta}V \\
\Sigma_o & @<{q}<< & \Sigma
\end{CD}
$$
denote the diagram where $q$ is the composition of the minimal sequence of blow-ups such that the pull-back $\tilde Z$ is smooth. 
Let $\tilde {S}_k$ be the
$\CC(S_k)$-normalization of $\Sigma$. The variety $\tilde {S}_k$ is a $D_{2k}$-cover 
of $\Sigma$ and we denote the cover morphism by $\tilde {\pi}_k$. Summing up, we 
obtain the following commutative diagram:
$$
\begin{CD}
S_k & @<<< & \tilde {S}_k \\
@VVV &  & @VVV \\
Z & @<{\mu}<< & \tilde Z \\
@V{\delta}VV&  & @VV{\tilde\delta}V \\
\Sigma_o & @<{q}<< & \Sigma.
\end{CD}
$$

Note that 
$$
\Delta_{\tilde \delta} =  q^{-1}C_1 + 
\text{ Some irreducible components of the exceptional set of $q$,}
$$
$$
\Delta_{\beta_2(\tilde{\pi}_k)}  =  \tilde \delta^{-1}(q^{-1}C_2) + 
\text{Some irreducible components of the exceptional set of $\mu$,} 
$$
where $\bullet^{-1}$ denote proper transforms.

By Corollary \ref{cor:branch}, $\tilde \delta^{-1}(q^{-1}C_2)$ is of the form $C_2^+ + C_2^-$,
$\sigma_{\tilde \delta}^*(C_2^+) = C_2^-$. Since $\pi_k$ is branched at $2C_1 + kC_2$, 
by Proposition \ref{prop:necdi}, for all $k$ as in the statement
there exists a line bundle $L_k$ such that
\begin{equation}
\label{ast}
C_2^+ - C_2^- + R_k - \sigma_{\tilde \delta}^*R_k \sim kL_k
\end{equation}
where $\Supp(R_k\cup\sigma_{\tilde \delta}^*R_k)$ is contained in the exceptional 
set of $\mu$.
The subgroup $T$ of $\NS(\tilde Z)$ generated by the irreducible components of the 
exceptional divisors of $\mu$ is a negative definite sublattice in $\NS(\tilde Z)$. 
Let us consider the relation $\eqref{ast}$ in $\NS(\tilde Z)/T$. Then we have 
\begin{equation}
\label{ast1}
C_2^+ - C_2^- \equiv kL_k \bmod T.
\end{equation}
Since $\NS(\tilde Z)/T$ is a finitely generated Abelian group, 
the hypothesis implies that $C_2^+ - C_2^-$ is a torsion element
of $\NS(\tilde Z)/T$; we can apply Remark~\ref{rem-num} since
$C_2^\pm$ is orthogonal to $T$.
Hence there exists a certain $\ell \in \NN$ such that 
$\ell(C_2^+ - C_2^-)\in T$. Put
$$
\ell(C_2^+ - C_2^-)= \sum_i c_i\Theta_i,
$$
where $\Theta_i$'s denote the irreducible components of the exceptional divisor of $\mu$. 
Since $C_2$ does not pass through the singularities of $C_1$
then $\Theta_i\cdot C_2^{\pm} = 0$ for all $i$. 
Hence $\ell(C_2^+ - C_2^-) = 0$ and as $T$ is a free $\ZZ$-module generated by $\Theta_i$'s, then
$C_2^+ = C_2^-$ in $\NS(\tilde Z)$. 
Since $\tilde Z$ is simply connected, $\Pic(\tilde Z) = \NS(\tilde Z)$. 
This implies $C_2^+ - C_2^- \sim 0$. Hence by
Corollary \ref{cor:con_di-1}, our statement follows.
\end{proof}

%\input inf_dihed_sec3.tex
%inf_dihed_sec3.tex

\section{Proof of Theorem \ref{thm:main}}

Let $\delta : Z \to \PP^2$ be a double cover branched at $C_1$, and let 
$\mu : \tilde Z \to Z$ be its canonical resolution. Since $C_1$ has at most
simple singularities, $\tilde Z$ is simply connected by Lemma~\ref{lem:P2}. Hence 
the first half of Theorem~\ref{thm:main} follows from Proposition~\ref{prop:key}.

We now go on to the second half. Since $C_2^+ - C_2^- \sim 0$
we can construct a family
$\pi_n : S_n \to \PP^2$, $n\in\NN$, of $D_{2n}$-covers branched at
$2C_1 + nC_2$ such that: 

\begin{enumerate}[(i)]
\smallbreak\item $\beta_1(\pi_n) = \delta : Z \to \PP^2$ is a double cover branched 
as above, and
\smallbreak\item $\CC(S_n) = \CC(Z)(\sqrt[n]{\varphi})$, where $\varphi$ is a rational
function with $(\varphi) = C_2^+ - C_2^-$.
\end{enumerate}

Note that $\CC(Z) = \CC(\PP^2)(\varphi)$ and 
$\CC(S_n) = \CC(\PP^2)(\sqrt[n]{\varphi})$. 
Let $[U: V: W]$ be a generic system of homogeneous coordinates of $\PP^2$. Assume that $C_1$ and $C_2$ are given by the equations:
\begin{align*}
C_1 : F_1(U, V, W) & =  0 \\
C_2 : F_2(U, V, W) & =  0,
\end{align*}
and let $L_{\infty} = \{ W = 0\}$. We may assume that both $C_1$ and $C_2$ meet
$L_{\infty}$ transversely. Put $u := U/W$, $v := V/W$, $f_1(u, v):= F_1(u, v, 1)$, and
$f_2(u, v) := F_2(u, v, 1)$. Let
$$
\scU := \PP^2 \setminus L_{\infty} = \Spec (\CC[u, v])\cong\CC^2,
$$
and
$$
\scV := {\delta}^{-1}(\scU) = \Spec(\CC[u, v, \zeta]), \quad \zeta^2 = f_1.
$$
Since $C_2 \cap \Sing(C_1) = \emptyset$, both $C_2^+$ and $C_2^-$ are Cartier divisors.
Let $L_{C_2^+}$ and $L_{C_2^-}$ be line bundles corresponding to 
$C_2^+$ and $C_2^-$. Let us consider $\scV_1$ an affine open set on $\scV$ invariant
under the involution $\sigma_\delta$ and such that
both ${L_{C_2^+}}|_{\scV_1}$ and $L_{C_2^-}|_{\scV_1}$ trivialize.

Then $\scV_1 = \Spec (R)$, for some ring $R$ satisfying
$\CC[u, v, \zeta] \subset R \subset \CC(u, v, \zeta)$.
Then on $\scV_1$ we may assume that $\scV_1\cap C_2^+$ and $\scV_1\cap C_2^-$ are
given by the local equations:
\begin{eqnarray*}
\scV_1\cap C_2^+ : \frac {g_1}a + \frac {g_2}a \zeta & = & 0 \\
\scV_1\cap C_2^- : \frac {g_1}a - \frac {g_2}a \zeta & = & 0,
\end{eqnarray*}
where $g_1, g_2, a \in \CC[u, v]$ and $\gcd(g_1, g_2) = 1$ (since $C_2^{\pm}$ 
are irreducible).
Hence the function $g_1^2 - g_2^2\zeta^2 = g_1^2 - g_2^2 f_1$ vanishes on $C_2\cap \delta(\scV_1)$,
that is $f_2 \mid (g_1^2 - g_2^2f_1)$. Thus we have 
$$
f_2 h = g_1^2 - g_2^2f_1 
$$
for some $h \in \CC[u, v]$. Homogenizing one obtains
$$
F_2(U, V, W)H(U, V, W)W^k = G_1^2 - G_2^2F_1,
$$
where $H(U,V,W) = W^{\deg h}h(U/W, V/W)$. Our statement is hence a 
consequence of the following claim.

\par\medskip

\begin{claim}
Both $k$ and $\deg H$ are $0$.
\end{claim}

\begin{proof}[Proof of Claim] Put $\theta_n = \sqrt[n]{\varphi}$ and
consider the $D_{2n}$-equivariant rational map
\[
\mu_n : S_n \dasharrow \PP^1
\]
given by $\theta_n$. By our construction, $\CC(S_n) = \CC(\PP^2)(\theta_n)$ for
all $n$. Hence  $S_n$ is regarded as a rational pull-back of 
$\PP^1 \to \PP^1, t \mapsto s:= t^n +t^{-n}$ by 
$$
s = \varphi + \frac{1}{\varphi} = \frac {2(G_1^2 + G_2^2 F_1)}{G_1^2 - G_2^2F_1}.
$$
Also there is no cancellation between the denominator and the numerator in
the last term in the above equality. Therefore if either $k \neq 0$ or 
$\deg H > 0$, then either $L_{\infty}$ or the curve given by $H(U, V, W) = 0$ 
is contained in the branch locus of $\pi_n$ for some $n$. 
This contradicts our assumption.
\end{proof}
%inf_dihed_sec4.tex
\section{Pencils and fundamental groups}

Let $C$ be a complex projective plane curve.
In this section we intend to exhibit the connection between the existence of
pencils of curves related to $C$ and the fundamental group
of its complement $X_C:=\PP^2\setminus C$ from a topological point of view. We will apply it to curves satisfying the statement
of Theorem~\ref{thm:main}.

\begin{defin} Let $D$ be a compact algebraic curve, let $p_1,\dots,p_r,q_1,\dots,q_s\in D$ be distinct points and
let $n_1,\dots,n_r\in\ZZ_{\geq 2}$.
An \emph{orbifold} $D^{q_1,\dots,q_s}_{(p_1,n_1),\dots,(p_r,n_r)}$
is a punctured curve $D\setminus\{q_1,\dots,q_s\}$
where the points $p_i$ are weighted with the integers $n_i$, $i=1,\dots,r$. For the sake of simplicity sometimes it will be
denoted by $D^{q_1,\dots,q_s}_{n_1,\dots,n_r}$.
\end{defin}

We may think that the charts around the points $p_i$ are obtained
as the quotient of disks in $\CC$  by the action of the $n_i$-roots
of unity. This justifies the following definition.

\begin{defin} The {orbifold-fundamental group} $\pior(D^{q_1,\dots,q_s}_{(p_1,n_1),\dots,(p_r,n_r)};*)$,
$*\in \check D:=D\setminus\{p_1,\dots,p_r,q_1,\dots,q_s\}$
is defined as the quotient of $\pi_1(\check D;*)$
by the normal subgroup generated by
$\mu_i^{n_i}$, $i=1,\dots,r$, where $\mu_i$
is a meridian of $p_i$.
\end{defin}

\begin{exmples}
We will fix $D=\PP^1$.
\begin{enumerate}
\item $\pior((\PP^1)_{2,2,n};*)$ is the dihedral group $D_{2 n}$.
\item $\TT_{p,q,r}:=\pior((\PP^1)_{p,q,r};*)$ is the corresponding triangle group.
\item $\FF_{n_1,\dots,n_r}:=\pior((\PP^1)^{\infty}_{n_1,\dots,n_r};*)$ is the free
product $\ZZ/n_1*\dots*\ZZ/n_r$.
\end{enumerate}
\end{exmples}

Let us fix now a connected smooth projective surface $X$, a connected smooth projective curve $\Gamma$ and a non-constant rational map $\widetilde\rho:X\dasharrow\Gamma$. Let $C\subset X$ be a compact curve such that $\widetilde\rho$ is well defined on $X\setminus C$ and let
$A:=\Gamma\setminus\widetilde\rho(X\setminus C)$, which is a finite set of points. We denote by $\rho:X\setminus C\to\Gamma\setminus A$ the restriction of $\widetilde\rho$, which is
assumed to have connected fibers. 

Let $p\in\Gamma\setminus A$; we consider the divisor $\rho^*(p)$, which is the restriction of $\widetilde\rho^*(p)$ to $X\setminus C$. For each $p$ we denote $n_p$ the $\gcd$ of the multiplicities of $\rho^*(p)$.
We consider the orbifold $\Gamma_\rho:=\Gamma^A_{\{(p,n_p)\mid n_p>1\}}$.
Fix $q\in\Gamma\setminus A$ such that $n_q=1$ and $*\in\rho^{-1}(q)$.

\begin{prop}\label{haz-orb} The mapping $\rho$ induces a natural epimorphism $\rho_*:\pi_1(X\setminus C;*)\twoheadrightarrow\pior(\Gamma_\rho;q)$.
\end{prop}

\begin{proof} Let us denote 
$\widetilde C:=C\cup\bigcup_{n_p>1}\widetilde\rho^*(p)$ and 
$\Gamma_1:=\Gamma\setminus(A\cup\{(p,n_p)\mid n_p>1\})$.
The rational map $\widetilde\rho$ induces
a well-defined surjective morphism 
$\rho_1:X\setminus\widetilde C\to\Gamma_1$.
It is a standard fact that $\rho$ induces an epimorphism 
$$
\pi_1(X\setminus \widetilde C;*)\twoheadrightarrow\pi_1(\Gamma_1;q).
$$
Recall that $\pi_1(X\setminus C;*)$ is the quotient
of $\pi_1(X\setminus \widetilde C;*)$ by the subgroup generated by the components of $\widetilde C$ not in $C$. The condition on the $\gcd$ of multiplicities guarantees the following commutative diagram which gives the result:
$$
\begin{CD}
\pi_1(X\setminus \widetilde C;*)&@>>> &\pi_1(\Gamma_1;q)\\
@VVV&&@VVV\\
\pi_1(X\setminus C;*)&@>>> &\pior(\Gamma_\rho;q).
\end{CD}
$$
Let us note that a meridian of a component of $\widetilde C$ not in $C$
is sent by $\rho$ to the power of a meridian $\mu_i$; the power
is a multiple of $n_i$.
\end{proof}

We say that a pencil $\scP:=\{F_p\}_{p\in\PP^1}$ {\emph{contains}} $C$ if each irreducible component of
$C$ is contained in a member of $\scP$. Let $A\subset\PP^1$ the subset
of $p\in\PP^1$ such that $F_p^{\text{red}}\subset C$. Let $n_p$ denote the $\gcd$ of 
the multiplicities of the components in $F_p$ not contained in $C$. We define the set
$B=\{p\in \PP^1\setminus A \mid n_p>1 \}\subset\PP^1$. Let us assume that $\#A=n$ and
$B:=\{p_1,\dots,p_r\}$, $n_i:=n_{p_i}$.

\begin{cor}
\label{thm:pencil}
There is a surjection from $\pi_1(X_C)$ onto 
$$
\FF_{n;(n_1,\dots,n_r)}:=\langle 
x_1,...,x_n,y_1,\dots,y_r : \prod_{j=1}^n x_j\cdot\prod_{i=1}^r y_i=y_1^{n_1}=...=y_r^{n_r}=1
\rangle
$$
\end{cor}

\begin{rem}
\label{rem:quo}
If $n'_i|n_i$, then $\FF_{n;(n_1,...,n_r)}$ surjects onto $\FF_{n;(n'_1,...,n'_r)}$. 
Any $n'_i$ equal~$1$ will be dropped. By doing so, we only add some ambiguity about
the surjection, but this is not relevant for our purposes.
\end{rem}

\begin{exmple}\label{exam-zar1}
Let $C_6$ be a Zariski sextic, that is, of equation $D^3_2+D^2_3=0$, where $D_i$ is
a homogeneous polynomial in $\CC[x,y,z]$ of degree $i$. The pencil generated by 
$D_2^3$ and $D_3^2$ has at least these two as special fibers. According to the notation 
of Corollary~\ref{thm:pencil} we have that
$\pi_1(X_{C_6})$ surjects, onto a group 
$\FF_{1;(2,3,n_3,\dots,n_r)}$ and
therefore (Remark~\ref{rem:quo}) onto 
$\FF_{1;(2,3)}=\ZZ_2 * \ZZ_3$.
Zariski proved in \cite{zariski1} that this is an isomorphism
for generic choices.
\end{exmple}

\begin{proof}[Proof of Corollary~\ref{cor:main}]
By Theorem~\ref{thm:main} the pencil generated by $G_1^2$ and $G_2^2F_1$ contains $F_2$,
therefore, using Corollary~\ref{thm:pencil}, there exists a surjection from 
$\pi_1(\PP^2\setminus (C_1\cup C_2))$ onto $\FF_{1;(2,2,n_3,...,n_r)}$, and hence,
by Remark~\ref{rem:quo}, there exists a surjection onto 
$\FF_{1;(2,2)}=\ZZ_2*\ZZ_2$.
\end{proof}

%section examples
\section{Examples}

\begin{exmple}
Let us suppose that there exists a pencil with three
fibers $2A_1+B_1$, $2A_2+B_2$, $n A_3+B_3$.
Then the fundamental group of the complement of
$B_1\cup B_2\cup B_3$ surjects onto $\TT_{2,2,n}=D_{2n}$.
The simplest example is the tricuspidal quartic.
Zariski proved in \cite{zariski1} that it lives
in a pencil as in Example~\ref{exam-zar1} if we add
the double of the bitangent line. Then we have
a surjection onto $D_6$ which is not an isomorphism
since it is also proved in \cite{zariski1} that its
fundamental group has order~$12$.
\end{exmple}

\begin{exmple}
Let $C$ be a smooth conic and $L_1$, $L_2$, $L_3$ tangent lines at three different
points $P_1,P_2$ and $P_3$ of $C$. The pencil $\scP$ generated by $C$ and $L_1+L_2$ 
contains as a special fiber $2L$, where $L$ is the line passing through $P_1$ and $P_2$.
Let $f_n$ be the cover map $f_n:\PP^2 \to \PP^2$, $f_n(L_1,L_2,L_3):=[L_1^n:L_2^n:L_3^n]$.
The pull-back $f^*\scP$ of the pencil $\scP$ is generated by $f_n^*C$ and 
$f_n^*L_1+f_n^*L_2=n(L_1+L_2)$ and contains the curve $2f_n^*L$. A description of the curve
$f_n^*C$ and a presentation of its fundamental group $\pi_1(\PP^2\setminus f_n^*C)$ can be
found in~\cite{ji}. By Corollary~\ref{thm:pencil}, $\pi_1(\PP^2\setminus f_n^*C)$ has a 
surjection onto $\FF_{1;(2,n)}=\ZZ_2*\ZZ_n$.
\end{exmple}

\begin{exmple}
In \cite{acct:01} we have studied curves having two irreducible components: a quartic $C_1$
having two singular points of types $\bbA_3$ and $\bbA_1$ and a smooth conic $C_2$
such that its intersection with $C_1$ produces a singular point of type $\bbA_{15}$.
Let us drop the $\bbA_1$ point. Then it is easily seen that the moduli space of such curves
has three connected components. Let us describe two of them:
\begin{itemize}
\item The tangent line $T$ at $\bbA_{15}$ passes through $\bbA_3$. In this case
there is a pencil of quartics containing $C_1$ and $4T$ such that another
element of the pencil is $C_2+2 L$, where $L$ is the tangent line at $\bbA_3$.
By Corollary~\ref{thm:pencil}, $\pi_1(\PP^2\setminus(C_1\cup C_2))$ has a 
surjection onto $\FF_{1;(2,4)}=\ZZ_2*\ZZ_4$.

\item There exists a smooth conic $Q$ having four infinitely near points in common with $\bbA_{15}$ and tangent at $\bbA_3$. In this case
there is again a pencil of quartics containing $C_1$, $2Q$ and $C_2+2 L$.
Therefore, $\pi_1(\PP^2\setminus(C_1\cup C_2))$ has a 
surjection onto $\FF_{1;(2,2)}=\ZZ_2*\ZZ_2$.
\end{itemize}
\end{exmple}

\begin{exmple}
Let us consider the family of curves of type $I$ described in~\cite{act}. When
$D$ is rational, they satisfy the conditions of Theorem~\ref{thm:main} and therefore
there is a surjection $\pi_1(\PP^2\setminus (D\cup L_1\cup L_2))\surj  \ZZ_2*\ZZ_2$ 
(Corollary~\ref{cor:main}). Note that there is yet another pencil that produces
a surjection $\pi_1(\PP^2\setminus (D\cup L_1\cup L_2))\surj  \ZZ_2*\ZZ_2$. Consider
the most general case, that is:
\begin{enumerate}
\item $C$ a rational arrangement of degree $2k+1$ with an ordinary multiple point $P$ 
of multiplicity $2k-1$ and at least $2k-1$ nodes $Q_1,...,Q_{2k-1}$.
\item $L_i$ a line tangent to $C$ at a point $P_i$, $i=1,2$.
\item $D_i$ a curve of degree $k$ with an ordinary multiple point at $P$ of 
multiplicity $k-1$, passing through $P_i,Q_1,...,Q_{2k-1}$.
\end{enumerate}
The pencil generated by $L_1+2D_2$ and $L_2+2D_1$ contains $C$.
Using a third line $L_3$ and the cover $f_n:\PP^2\to \PP^2$ described above,
one obtains curves $f_n^*C$ whose fundamental group $\pi_1(\PP^2\setminus f_n^*C)$
surjects onto $\ZZ_2*\ZZ_2$ (for $n$ even) and such that 
$\pi_1(\PP^2\setminus (f_n^*C\cup f_n^*D_1 \cup f_n^*D_2))$ surjects onto 
$\ZZ_n*\ZZ_n$.
\end{exmple}


\begin{thebibliography}{99}
\providecommand{\bysame}{\leavevmode\hbox to3em{\hrulefill}\thinspace}

\bibitem{acct:01}
E.~Artal, J.~Carmona, J.I. Cogolludo, and H.~Tokunaga, \emph{Sextics with
  singular points in special position}, J. Knot Theory Ramifications
  \textbf{10} (2001), no.~4, 547--578.

\bibitem{act} E. Artal, J.I. Cogolludo-Agust{\'\i}n and H. Tokunaga: 
\emph{Nodal degenerations of plane curves and Galois covers}, 
\rm preprint, 2004.

\bibitem{brieskorn1} E. Brieskorn: \emph{\"Uber die Aufl\"osung gewisser Singularit\"aten von
holomorpher Abbildungen}, \rm Math. Ann. \textbf{166} (1966), 76--102.

\bibitem{brieskorn2} \bysame: \emph{\"Uber die Aufl\"osung der rationalen Singularit\"aten 
holomorpher Abbildungen}, \rm Math. Ann. \textbf{178} (1968), 255--270.

\bibitem{ji} J.I. Cogolludo-Agust{\'\i}n: 
\emph{Fundamental group for some cuspidal curves}, \rm Bull. London Math. Soc. 
\textbf{31} (1999), no.~2, 136--142.

\bibitem{horikawa} E. Horikawa: \emph{On deformation of quintic surfaces},
\rm Invent. Math. {\bf 31} (1975), 43--85.

\bibitem{tokunaga1} H. Tokunaga: \emph{On dihedral Galois covers}, Canadian
J. of Math. {\bf 46} (1994), 1299--1317. 

\bibitem{tokunaga2}\bysame: \emph{Dihedral covers of algebraic surfaces and
its application}, Trans. Amer. Math.Soc. \textbf{352} (2000), 4007--4017.

\bibitem{zariski1} O. Zariski:  \emph{On the problem of existence of algebraic 
functions of two variables possessing a given branch curve}, Amer. J. Math.
\textbf{51} (1929), 305--328.

\bibitem{zariski4} \bysame: \emph{On the purity of the branch locus of algebraic functions}, 
Proc. Nat. Acad. USA \textbf{44} (1958), 791--796.

\end{thebibliography}
\end{document}